\documentclass{amsart}
\usepackage{amsfonts}
\usepackage{amssymb}
\usepackage{amsmath}

\newtheorem{theorem}{Theorem}
\theoremstyle{plain}

\newtheorem{corollary}[theorem]{Corollary}

\newtheorem{lemma}[theorem]{Lemma}

\newtheorem{proposition}[theorem]{Proposition}
\newtheorem{remark}{Remark}

\begin{document}
\title[] { On the $q$-extensions of the Bernoulli and Euler numbers, related identities and Lerch zeta function }
\author{Taekyun Kim  }
\address{Taekyun Kim. Division of General Education-Mathematics, \\
Kwangwoon University, Seoul 139-701, Republic of Korea  \\}
\email{tkkim@kw.ac.kr}
\author{Young-Hee Kim}
\address{Young-Hee Kim. Division of General Education-Mathematics, \\
Kwangwoon University, Seoul 139-701, Republic of Korea  \\}
\email{yhkim@kw.ac.kr}
\author{Kyung-Won Hwang}
\address{Kyung-Won Hwang. Department of General education,\\
Kookmin university, 861-1 Seongbukgu Seoul 136-702, Republic of
Korea \\} \email{khwang7@kookmin.ac.kr}

\maketitle

{\footnotesize {\bf Abstract} \hspace{1mm} {Recently,
$\lambda$-Bernoulli and $\lambda$-Euler numbers are studied in [5,
10]. The purpose of this paper is to present a systematic study of
some families of the $q$-extensions of the $\lambda$-Bernoulli and
the $\lambda$-Euler numbers by using the bosonic $p$-adic
$q$-integral and the fermionic $p$-adic $q$-integral. The
investigation of these $\lambda$-$q$-Bernoulli and
$\lambda$-$q$-Euler numbers leads to interesting identities related
to these objects. The results of the present paper cover earlier
results concerning $q$-Bernoulli and $q$-Euler numbers. By using
derivative operator to the generating functions of
$\lambda$-$q$-Bernoulli and $\lambda$-$q$-Euler numbers, we give the
$q$-extensions of Lerch zeta function.  }}



\section{Introduction, Definitions and Notations }

Throughout this paper, the symbols $\mathbb{Z}_p, \mathbb{Q}_p,
\mathbb{C}$ and $\mathbb{C}_p$ denote the ring of $p$-adic rational
integers, the field of $p$-adic rational numbers, the complex number
field and the completion of algebraic closure of $\mathbb{Q}_p$,
respectively. Let $\mathbb{N}$ be the set of natural numbers.

The symbol $q$ can be treated as a complex number, $q \in
\mathbb{C}$, or as a $p$-adic number, $q \in \mathbb{C}_p$. If $q
\in \mathbb{C}$, then we always assume that $|q|< 1.$ If $q \in
\mathbb{C}_p$, then we usually assume that $|1-q|_p< 1$. Here
$|\cdot|_p$ stands for the $p$-adic absolute value in $\mathbb{C}_p$
with $|p|_p = \frac{1}{p}$. The $q$-basic natural numbers are
defined by $[n]_q = \frac{1-q^n}{1-q}= 1+q+q^2+ \cdots + q^{n-1}~ (
n \in \mathbb{N})$ and $[n]_{-q} = \frac{1-(-q)^n }{1+q}$. In this
paper, we use the notation $$[x]_q = \frac{1-q^x}{1-q} \quad
\text{and} \quad [x]_{-q} = \frac{1-(-q)^x}{ 1+q}, \text{ see
[1-19]}.$$ Hence $ \underset{q \rightarrow 1} {\lim} [x]_q= x$ for
any $x$
 with $|x|_p \leq 1$ in the present $p$-adic case.

For $x \in \mathbb{Z}_p$, we say that $g$ is a uniformly
differentiable function at a point $a \in \mathbb{Z}_p,$ and write
$g \in UD (\mathbb{Z}_p),$ the set of uniformly differentiable
function, if the difference quotients
$$F_g (x,y)=\frac{g(y)-g(x)}{y-x}$$ have a limit $l=g'(a)$ as $(x,
y)\rightarrow (a,a)$. For $f \in UD(\mathbb{Z}_p )$, the $q$-
deformed bosonic $p$-adic integral is defined as
\begin{eqnarray}
I_q(f)= \int_{\mathbb{Z}_p} f(x)d\mu_q (x)=\underset{N
\rightarrow\infty} {\lim} \underset {x=0}{\overset{p^{N}-1}{\sum}}
f(x)\frac{q^x}{[p^{N}]_q} \label{Iq}, \text{ see
[1-19]},\end{eqnarray} and the $q$-deformed fermonic $p$-adic
integral is defined by
\begin{eqnarray*}
I_{-q} (f)= \int_{\mathbb{Z}_p} f(x)d\mu_{-q} (x)=\underset{N
\rightarrow\infty} {\lim} \underset {x=0}{\overset{p^{N}-1}{\sum}}
f(x)\frac{(-q)^x}{[p^{N}]_{-q}}, \quad  (\text{see \, [1-19]}).
\end{eqnarray*}
For $n \in \mathbb{N}$, let $f_n (x)= f(x+n).$ Then
\begin{eqnarray} q^n I_{-q}(f_n)= (-1)^n I_{-q}(f)+ [2]_q \underset
{l=0}{\overset{n-1}{\sum}} (-1)^{n-1-l} q^l f(l).
\end{eqnarray}

The classical Bernoulli polynomials $B_n(x)$ and the Euler
polynomials $E_n(x)$ are defined as
\begin{eqnarray}
\frac{t}{e^t-1}e^{xt}= \underset{x=0}{\overset{\infty}{\sum}}
B_n(x)\frac{t^n}{n!} \quad \text{and} \quad \frac{2e^{xt}}{e^t +1}=
\underset{x=0}{\overset{\infty}{\sum}} E_n (x)\frac{t^n}{n!}.
\end{eqnarray}
The Bernoulli numbers $B_n$ and the Euler numbers $E_n$ are defined
as $B_n = B_n(0)$ and $E_n= E_n(0)$, (see [1-19]).

From (\ref{Iq}), we note that
\begin{eqnarray}
q I_{q}(f_1)=  I_{q}(f)+ (q-1)f(0)+ \frac{q-1}{\log q}f'(0),
\label{Iq1}
\end{eqnarray}
for $f_1 (x)=f(x+1)$. By (\ref{Iq1}), we see that $I_1 (f_1)= I_1
(f)+f'(0),$ (see [7]).

Let $u$ be algebraic in $\mathbb{C}_p$ (or $\mathbb{C}$). Then the
Frobenius-Euler polynomials are defined as
\begin{eqnarray}
\frac{1-u}{e^t-u} \, e^{xt}= \underset{n=0}{\overset{\infty}{\sum}}
H_n(u, x)\frac{t^n}{n!}, \quad (\text{see} \,\, [5]).
\end{eqnarray}
In case $x=0$, $H_n (u,0)= H_n(u)$, which are called the Frobenius
Euler numbers.

Let $C_{p^n}$ be the cyclic group consisting of all $p^n$-th roots
of unity in $\mathbb{C}_p$ for any $n \geq 0$ and $T_p$ be the
direct limit of $C_{p^n}$ with respect to the natural morphisms,
hence $T_p$ is the union of all $C_{p^n}$ with discrete topology.

For $\lambda \in T_p$ with $\lambda\neq 1$, if we use (\ref{Iq1}),
then we have
\begin{eqnarray}
\int_{\mathbb{Z}_{p}} e^{tx} \lambda ^x d\mu_1 (x) =
\frac{t}{\lambda e^t -1}. \label{elam}
\end{eqnarray}

From (\ref{elam}), the $\lambda-$Bernoulli numbers are defined as
\begin{eqnarray}
\frac{t}{\lambda e^t -1}= e^{B(\lambda)t}=
\underset{n=0}{\overset{\infty}{\sum}} B_n (\lambda)\frac{t^n}{n!},
\quad (\text{see} \,\, [5])
\end{eqnarray}
with the usual convention of replacing $B^i (\lambda)$ by
$B_i(\lambda)$. Thus, $B_k(\lambda)$ can be determined inductively
by
\begin{eqnarray}
\lambda (B(\lambda)+1)^k - B_k(\lambda) =
\begin{cases}1, &\text{if $k = 1$,}\\
             0, &\text{if $k > 1$,} \quad(\text{see [5]}).
\end{cases}
\end{eqnarray}

By the definition of the Frobenius-Euler numbers, we see that
\begin{eqnarray}
\frac{t}{\lambda e^t-1}=\underset{m=0}{\overset{\infty}{\sum}}
\frac{1}{(m+1) !} \cdot \frac{(m+1)H_m(\lambda ^{-1})}{\lambda-1}
t^{m+1}, \quad\text{( see [7])}.
\end{eqnarray}
For $m\geq 1$ and $\lambda \neq 1$, we have
\begin{eqnarray}
B_m(\lambda) = \int _ {\mathbb{Z}_p} x^m \lambda ^x d\mu _1 (x)=
\frac{m}{\lambda -1} H_{m-1} (\lambda ^{-1}),\quad (\text{see [5]}).
\end{eqnarray}

We can also easily see that $\int _ {\mathbb{Z}_p} \lambda^x d\mu _1
(x)=0$ and
\begin{eqnarray*}
e^t x =\underset{m \rightarrow\infty} {\lim} \underset {\lambda \in
C_{p^m}}{\sum} \frac{t\lambda^x}{\lambda e^t-1}
=\underset{n=0}{\overset{\infty}{\sum}} \frac{t^n}{n!} \underset{m
\rightarrow\infty} {\lim} \underset {\lambda \in C_{p^m}}{\sum}\int
_ {\mathbb{Z}_p} x^m \lambda ^x d\mu _1 (x)\lambda ^x .
\end{eqnarray*}
Consequently, we have
\begin{eqnarray*}
x^n &=& B_n(1)+ \underset {\lambda \neq 1}{\underset {\lambda \in
T_{p}}{\sum}} \frac{1}{\lambda-1} H_{n-1} (\lambda ^{-1})
\lambda ^x\\
&=& B_n(1)+ \underset {\lambda \neq 1}{\underset {\lambda \in
T_{p}}{\sum}} \frac{B_n (\lambda)}{n} \lambda ^x, \quad\text{( see
[5])}.
\end{eqnarray*}
From $(6)$ and $(8)$, we note that $$B_0 (\lambda)=0, \, B_1
(\lambda)= \frac{1}{\lambda -1}, \, B_2 (\lambda)= -
\frac{2\lambda}{(\lambda-1)^2}, \, \cdots.$$

The Genocchi numbers are defined by the generating function
$$\frac{2t}{e^t +1}= \underset{n=0}{\overset{\infty}{\sum}} G_n
\frac{t^n}{n!}.$$ These numbers satisfy the relation $G_0=0, G_1= 1,
G_3= G_5 = \cdots = G_{2k+1} = 0$, and the even coefficients are
$G_n= 2(1-2^n)B_n .$

For $\lambda \in \mathbb{C}_p$ with $|\lambda|< 1$, by $(2)$, we
have
\begin{eqnarray}
\int _ {\mathbb{Z}_p} \lambda ^x e^{xt} d\mu _{-1} (x) =
\frac{2}{\lambda e^t +1}.
\end{eqnarray}

By $(11)$, we define the $\lambda$-Euler numbers as follows :
\begin{eqnarray}
\frac{2}{\lambda e^t +1}= \underset{n=0}{\overset{\infty}{\sum}}
\frac{E_n (\lambda)}{n!}t^n, \quad (\text{see} \,\, [7, 9, 10]).
\end{eqnarray}
Note that $E_n (\lambda)=\frac{2}{\lambda+1}H_n(- \lambda ^{-1}).$

From $(12)$, we can easily derive
\begin{eqnarray} \int_{\mathbb{Z}_p} x^n \lambda
^x  d\mu _{-1} (x)=E_n (\lambda)= \frac{2}{\lambda+1}H_n(-\lambda
^{-1}).
\end{eqnarray}

The $\lambda$-Genocchi numbers are also defined as
$$t\int_{\mathbb{Z}_p} x^n \lambda ^x  d\mu _{-1} (x)= \frac{2t}{\lambda e^t
+1}=\underset{n=0}{\overset{\infty}{\sum}} G_n (x)\frac{t^n}{n!}. $$
Thus, we have $G_0 (\lambda)=0, \, G_1
(\lambda)=\frac{2}{\lambda+1}, \, \cdots, \, E_n (\lambda)=
\frac{G_{n+1}(\lambda)}{n+1}.$

In this paper, we study the $q$-extension of $\lambda$-Bernoulli
number and $\lambda$-Euler numbers related to Lerch zeta function.
The purpose of this paper is to present a systematic study of some
families of the $q$-extension of the $\lambda$-Bernoulli and
$\lambda$-Euler numbers by using the bosonic $p$-adic $q$-integral
and the ferminionic $p$-adic $q$-integral. The investigation of
these $\lambda$-$q$-Bernoulli and $\lambda$-$q$-Euler numbers leads
to interesting identities related to these objects. The results of
the present paper cover earlier results concerning $q$-Bernoulli and
$q$-Euler numbers. By using derivative operator to the generating
functions of $\lambda$-$q$-Bernoulli and $\lambda$-$q$-Euler
numbers, we can give the $q$-extension of Lerch zeta function.

\vskip 20pt

\section{$q$-extension of $\lambda$-Bernoulli numbers and polynomials}
\vskip 10pt

For $\lambda \in T_p $,  let us consider the $q$-extension of
$\lambda$-Bernoulli numbers as follows.
\begin{eqnarray}
\beta_{k,q} (\lambda)= {\int_{\Bbb Z_p}}\lambda^x
[x]_{q}^{k}d\mu_{q}(x). \label{lamB}
\end{eqnarray}
From (\ref{lamB}), we note that
\begin{eqnarray*}
\beta_{k,q} (\lambda)&=&\lim_{N \rightarrow \infty}
\frac{1}{[p^N]_{q}} \sum_{x=0}^{p^N -1} \lambda^x [x]_{q}^{k} q^x \\
&=&\lim_{N \rightarrow \infty} \frac{1}{[p^N]_{q}} \sum_{x=0}^{p^N
-1} (\lambda q)^x (\sum_{l=0}^k \binom{k}{l} (-1)^l q^{lx} )
\frac{1}{(1-q)^k} \\
&=&\frac{1-q}{(1-q)^k}\sum_{l=0}^k \binom{k}{l} (-1)^l
\frac{1-(\lambda q^{l+1})^{p^N}}{1-\lambda q^{l+1}}\\
&=& \frac{1}{(1-q)^{k-1}} \sum_{l=0}^k \binom{k}{l} (-1)^l
\frac{l+1}{1-\lambda q^{l+1}} \, .
\end{eqnarray*}
Therefore, we obtain the following theorem.

\begin{theorem} For $k \in {\Bbb N} \cup \{ 0 \}$ and $\lambda \in T_p \,$, we have
$$\beta_{k,q} (\lambda)=\frac{1}{(1-q)^{k-1}} \sum_{l=0}^k \binom{k}{l} (-1)^l
\frac{l+1}{1-\lambda q^{l+1}} \, .$$
\end{theorem}

\vskip 10pt

Let $F(t, \lambda : q)$ be the generating functions of $\beta_{n,q}
(\lambda)$ with
$$F(t, \lambda  : q)=\sum_{n=0}^\infty \beta_{n,q} (\lambda) \frac{t^n}{n !} \, .$$
Then we have
\begin{eqnarray}
F(t, \lambda  : q)&=&\sum_{n=0}^\infty \beta_{n,q} (\lambda)
\frac{t^n}{n !}={\int_{\Bbb Z_p}}\lambda^x e^{[x]_{q} t}
d\mu_{q}(x )\notag \\
&=&\sum_{n=0}^\infty {\int_{\Bbb Z_p}}\lambda^x
[x]_{q}^{n}d\mu_{q}(x) \frac{t^n}{n !} \notag \\
&=&\sum_{k=0}^\infty \{ \frac{1}{(1-q)^{k-1}} \sum_{l=0}^k
\binom{k}{l} (-1)^l (l+1)\sum_{m=0}^\infty \lambda^m q^{(l+1)m} \}
\frac{t^k}{k!} \label{gfbeta} \\
&=&\sum_{k=0}^\infty  \frac{1}{(1-q)^{k-1}} \sum_{m=0}^\infty
\lambda^m \sum_{l=0}^k \binom{k}{l} (-1)^l (l+1) q^{(l+1)m}
\frac{t^k}{k!} \notag \\
&=&\sum_{k=0}^\infty  \frac{1}{(1-q)^{k-1}} \sum_{m=0}^\infty q^m
\lambda^m \sum_{l=0}^k l \binom{k}{l} (-1)^l  q^{lm}
\frac{t^k}{k!} \notag \\
& &+\sum_{k=0}^\infty  \frac{1}{(1-q)^{k-1}} \sum_{m=0}^\infty q^m
\lambda^m \sum_{l=0}^k \binom{k}{l} (-1)^l  q^{lm} \frac{t^k}{k!}.
\notag
\end{eqnarray}
Since $l \binom{k}{l} = k \binom{k-1}{l-1} $, the first term of  the
last equation in (\ref{gfbeta}) equals
\begin{eqnarray}
& & \sum_{m=0}^\infty q^m \lambda^m \{ \sum_{k=0}^\infty
\frac{1}{(1-q)^{k-1}}\sum_{l=1}^k
\binom{k-1}{l-1} (-1)^l  q^{lm} \} \frac{t^k}{(k-1)!} \notag \\
& & \qquad =-\sum_{m=0}^\infty q^{2m} \lambda^m \{ \sum_{k=1}^\infty
\frac{1}{(1-q)^{k-1}}\sum_{l=0}^{k-1} \binom{k-1}{l} (-1)^l q^{lm}
\} \frac{t^k}{(k-1)!} \label{gfbeta1} \\
& & \qquad = -t \sum_{m=0}^\infty  q^{2m} \lambda^m
\sum_{k=0}^\infty [m]_q^k \frac{t^k}{k!}= -t \sum_{m=0}^\infty
q^{2m} \lambda^m e^{[m]_q t}. \notag
\end{eqnarray}
The second term of the last equation in (\ref{gfbeta}) equals
\begin{eqnarray}
& & \sum_{k=0}^\infty  \frac{1}{(1-q)^{k-1}} \sum_{m=0}^\infty q^m
\lambda^m (1-q^m )^k \frac{t^k}{k!}
\notag \\
& & \qquad = (1-q) \sum_{m=0}^\infty q^m \lambda^m
\sum_{k=0}^\infty [m]_q^k \frac{t^k}{k!} =(1-q)\sum_{m=0}^\infty q^m
\lambda^m e^{[m]_q t}. \label{gfbeta2}
\end{eqnarray}

From (\ref{gfbeta}), (\ref{gfbeta1}) and (\ref{gfbeta2}), we obtain
the following proposition.

\begin{proposition} Let $F(t, \lambda  : q)=\sum_{n=0}^\infty \beta_{n,q} (\lambda) \frac{t^n}{n !} \, .$
Then we have
$$F(t, \lambda  : q)= -t \sum_{m=0}^\infty
q^{2m} \lambda^m e^{[m]_q t}+(1-q)\sum_{m=0}^\infty q^m \lambda^m
e^{[m]_q t}.$$
\end{proposition}
Since $q^{2m} =q^m \{ [m]_q (q-1) +1 \} $, it follows that
\begin{eqnarray*} \beta_{k,q} (\lambda)&=& \frac{d^k F_q (t,
\lambda : q)}{(dt)^k
}|_{t=0} \\ &=& -k \sum_{m=0}^\infty q^{2m} \lambda^m [m]_q^{k-1}+ (1-q) \sum_{m=0}^\infty q^{m} \lambda^m [m]_q^{k} \\
&=& -k(q-1) \sum_{m=0}^\infty q^m \lambda^m [m]_q^k -
k\sum_{m=0}^\infty q^m \lambda^m [m]_q^{k-1} +(1-q)
\sum_{m=0}^\infty q^{m} \lambda^m [m]_q^{k}  \\
&=& (1-q)(k+1) \sum_{m=0}^\infty q^{m} \lambda^m [m]_q^{k} -
k\sum_{m=0}^\infty q^m \lambda^m [m]_q^{k-1} .
\end{eqnarray*}
Therefore, we obtain the following theorem.

\begin{theorem} For $k \in {\Bbb N} \cup \{ 0 \}$ and $\lambda \in T_p \,$, we have
$$\beta_{k,q} (\lambda)=(1-q)(k+1) \sum_{m=0}^\infty q^{m} \lambda^m [m]_q^{k} -
k\sum_{m=0}^\infty q^m \lambda^m [m]_q^{k-1} .$$
\end{theorem}

Now we consider another $q$-extension of $\lambda$-Bernoulli numbers
as follows.
\begin{eqnarray}
B_{n,q} (\lambda)= {\int_{\Bbb Z_p}} q^{-x}\lambda^x [x]_{q}^{n}
d\mu_{q}(x). \label{lamB1}
\end{eqnarray}
From (\ref{lamB1}), we can derive
\begin{eqnarray*}
B_{n,q} (\lambda)&=&{\int_{\Bbb Z_p}} q^{-x}\lambda^x  [x]_{q}^{n}
d\mu_{q}(x) \\
&=&\frac{1}{(1-q)^n}\sum_{l=0}^n \binom{n}{l} {\int_{\Bbb Z_p}}
q^{-x} (-1)^l \lambda^x q^{lx}d\mu_{q}(x)\\
&=& \frac{1}{(1-q)^{n-1}} \sum_{l=0}^n \binom{n}{l} (-1)^l
\frac{l}{1-\lambda q^{l}} \, .
\end{eqnarray*}
Thus, we obtain the following theorem.

\begin{theorem} For $n \in {\Bbb N} \cup \{ 0 \}$ and $\lambda \in T_p \,$, we have
$$B_{n,q} (\lambda)=\frac{1}{(1-q)^{n-1}} \sum_{l=0}^n \binom{n}{l} (-1)^l
\frac{l}{1-\lambda q^{l}} \, .$$
\end{theorem}

Let $F^*(t, \lambda : q)$ be the generating functions of $B_{n,q}
(\lambda)$ with
$$F^*(t, \lambda  : q)=\sum_{n=0}^\infty B_{n,q} (\lambda) \frac{t^n}{n !} \, .$$
Then we have
\begin{eqnarray*}
F^* (t, \lambda  : q)&=&\sum_{n=0}^\infty B_{n,q} (\lambda)
\frac{t^n}{n !}={\int_{\Bbb Z_p}}q^{-x} \lambda^x e^{[x]_{q} t}
d\mu_{q}(x ) \\
&=&\sum_{n=0}^\infty \{ {\int_{\Bbb Z_p}}q^{-x} \lambda^x
[x]_{q}^{n} d\mu_{q}(x) \} \frac{t^n}{n !}  \\
&=&\sum_{n=0}^\infty \{ \frac{1}{(1-q)^{n-1}} \sum_{l=0}^n
\binom{n}{l} (-1)^l \frac{l}{1-\lambda q^l } \} \frac{t^n}{n!} \label{gfB} \\
&=&\sum_{n=0}^\infty  \{ \frac{1}{(1-q)^{n-1}} \sum_{l=0}^n
\binom{n}{l}(-1)^l l \sum_{m=0}^\infty \lambda^m  q^{lm} \}
\frac{t^n}{n!} \\
&=&\sum_{m=0}^\infty  \lambda^m \{   \sum_{n=1}^\infty
\frac{n}{(1-q)^{n-1}}  \sum_{l=1}^n  \binom{n-1}{l-1} (-1)^l  q^{lm}
\} \frac{t^n}{n!}  \\
&=&-\sum_{m=0}^\infty  \lambda^m q^{m}  \sum_{n=1}^\infty
\frac{n}{(1-q)^{n-1}}  (1-q^m )^{n-1} \frac{t^n}{n!}  \\
&=&-\sum_{m=0}^\infty  \lambda^m q^{m}  \sum_{n=0}^\infty
\frac{(1-q^m )^{n}}{(1-q)^{n}}  \frac{t^{n+1}}{n!}  \\
&=&-t\sum_{m=0}^\infty  \lambda^m q^{m} e^{[m]_q t} .
\end{eqnarray*}
Therefore we obtain the following lemma.
\begin{lemma} Let $F^*(t, \lambda  : q)=\sum_{n=0}^\infty B_{n,q} (\lambda) \frac{t^n}{n !} \,
$.  Then we have
$$ F^*(t, \lambda  : q)=-t\sum_{m=0}^\infty  \lambda^m q^{m} e^{[m]_q t}.$$
\end{lemma}

We also have
$$B_{k,q} (\lambda)={ \frac{d^k F_q (t,
\lambda : q)}{(dt)^k }}|_{t=0}= - k\sum_{m=0}^\infty q^m \lambda^m
[m]_q^{k-1} .
$$
Therefore we obtain the following theorem.
\begin{theorem} For $k \in {\Bbb N} \cup \{ 0 \}$ and $\lambda \in T_p \,$, we have
$$B_{k,q} (\lambda)=-k \sum_{m=0}^\infty q^m \lambda^m
[m]_q^{k-1} .$$
\end{theorem}

\bigskip

\section{$q$-extension of $\lambda$-Euler numbers and polynomials }
\vskip 10pt

In this section, we assume that $p$ is an odd prime number and
$\lambda \in \mathbb{C}_p$ with $|1-\lambda|_p < 1$. By using the
fermionic $p$-adic $q$-integral on $\mathbb{Z}_p$, we consider the
$q$-extensions of $\lambda$-Euler numbers as follows.

For $n \in \mathbb{N}\cup \{ 0 \}$, we define the $q-$extension of
$\lambda$-Euler numbers as
\begin{eqnarray}
E_{n,q}(\lambda) =\int_{\mathbb{Z}_p} q^{-x} \lambda ^x [x]_q ^n
d\mu _{-q}(x). \label{lamE}
\end{eqnarray}

From (\ref{lamE}), we note that
\begin{eqnarray*}
E_{n,q}(\lambda)&=&\int _ {\mathbb{Z}_p} q^{-x}
\lambda ^x [x]_q ^n d\mu _{-q} (x)\\
&=&\underset{N \rightarrow\infty} {\lim} \frac{1+q}{\,\, 1+
q^{p^N}}\underset{x=0}{\overset{p^N -1}{\sum}} (-1)^x [x]_q ^n
\lambda^x \\
&=&
\frac{[2]_q}{2}\frac{1}{(1-q)^n}\underset{l=0}{\overset{n}{\sum}}{n
\choose l} (-1)^l \underset{N \rightarrow\infty} {\lim} \frac{1+
q^{p^N}\lambda^{p^N}}{1+q^l \lambda} \\
&=&
\frac{[2]_q}{2}\frac{1}{(1-q)^n}\underset{l=0}{\overset{n}{\sum}}{n
\choose l} (-1)^l \frac{2}{1+q^l \lambda}\\
&=&\frac{[2]_q}{(1-q)^n} \underset{l=0}{\overset{n}{\sum}}{n \choose
l} (-1)^l \frac{1}{1+q^l \lambda}.
\end{eqnarray*}

Therefore we obtain the following theorem.
\begin{theorem} For $n \in {\Bbb N} \cup \{ 0 \}$, we have
$$E_{n,q}(\lambda)= \frac{[2]_q}{(1-q)^n}
\underset{l=0}{\overset{n}{\sum}}{n \choose l} (-1)^l \frac{1}{1+q^l
\lambda}.$$
\end{theorem}

Let $g(t, \lambda: q)$ be the generating function of
$E_{n,q}(\lambda)$ with $$g(t, \lambda:
q)=\underset{n=0}{\overset{\infty}{\sum}}E_{n,q}(\lambda)\frac{t^n}{n!}.$$
Then we have
\begin{eqnarray*}
g (t, \lambda  : q)&=&\sum_{n=0}^\infty E_{n,q} (\lambda)
\frac{t^n}{n !}={\int_{\Bbb Z_p}}q^{-x} \lambda^x e^{[x]_{q} t}
d\mu_{-q}(x ) \\
&=&\sum_{n=0}^\infty \{ {\int_{\Bbb Z_p}}q^{-x} \lambda^x
[x]_{q}^{n} d\mu_{-q}(x) \} \frac{t^n}{n !}  \\
&=&[2]_q \sum_{n=0}^\infty \{ \frac{1}{(1-q)^{n}} \sum_{l=0}^n
\binom{n}{l} (-1)^l \frac{1}{1+\lambda q^l } \} \frac{t^n}{n!} \\
&=&[2]_q \sum_{n=0}^\infty  \frac{1}{(1-q)^{n}} \sum_{l=0}^n
\binom{n}{l}(-1)^l  \{ \sum_{m=0}^\infty (-1)^m \lambda^m  q^{lm} \}
\frac{t^n}{n!} \\
&=&[2]_q \sum_{m=0}^\infty  (-1)^{m} \lambda^m \sum_{n=0}^\infty
\frac{1}{(1-q)^{n}}  \sum_{l=0}^n  \binom{n}{l} (-1)^l  q^{lm} \frac{t^n}{n!}  \\
&=&[2]_q\sum_{m=0}^\infty (-1)^{m} \lambda^m  \sum_{n=0}^\infty
[m]_q^n  \frac{t^{n}}{n!}  \\
&=&[2]_q\sum_{m=0}^\infty (-1)^m \lambda^m e^{[m]_q t} .
\end{eqnarray*}

Thus, we have the following lemma.

\begin{lemma}
Let $g(t, \lambda:
q)=\underset{n=0}{\overset{\infty}{\sum}}E_{n,q}(\lambda)\frac{t^n}{n!}.$
Then we have
\begin{eqnarray}
g(t, \lambda: q)=[2]_q\sum_{m=0}^\infty (-1)^m \lambda^m e^{[m]_q
t}. \label{lamG}
\end{eqnarray}
\end{lemma}

By (\ref{lamG}), we can also consider the $\lambda$-$q$-Genocchi
numbers as follows.
\begin{eqnarray}
t {\int_{\Bbb Z_p}}q^{-x} \lambda^x e^{[x]_{q} t} d\mu_{-q}(x )
=[2]_q t \sum_{m=0}^\infty (-1)^m \lambda^m e^{[m]_q
t}=\sum_{n=0}^\infty G_{n,q}(\lambda)\frac{t^n}{n!}. \label{lamqG}
\end{eqnarray}

From (\ref{lamqG}), we note that $G_{0,q} (\lambda)=0$ and $$\int _
{\mathbb{Z}_p}q^{-x} \lambda^x  [x]_q ^n d\mu _{-q}
(x)=\frac{G_{n+1,q}(\lambda)}{n+1}.$$
Thus, we see that
$$E_{n,q}(\lambda)=\int_{\mathbb{Z}_p}q^{-x}
\lambda^x  [x]_q ^n d\mu _{-q} (x)=\frac{G_{n+1,q}(\lambda)}{n+1}.$$
Hence $$G_{n,q}(\lambda)=[2]_q
\frac{n}{(1-q)^{n-1}}\underset{l=0}{\overset{n-1}{\sum}}{n-1 \choose
l} (-1)^l \frac{1}{1+q^l \lambda},$$ where $n=1,  2, 3, \, \cdots$.
Indeed,
\begin{eqnarray*}
G_{1,q}(\lambda)&=&\frac{[2]_q}{1+\lambda},\\
G_{2,q}(\lambda)&=&\frac{2[2]_q}{1-q}\underset{l=0}{\overset{1}{\sum}}{1
\choose l} (-1)^l \frac{1}{1+q^l \lambda}=
\frac{[2]_q}{1-q}(\frac{2}{1+\lambda}-\frac{2}{1+q
\lambda})\\&=&-2[2]_q (\frac{\lambda}{(1+\lambda)(1+q\lambda)}).
\end{eqnarray*}

Now, we consider the $q$-extension of $\lambda$-Euler polynomials as
follows.
\begin{eqnarray}
E_{n,q}(\lambda,x)=\int _ {\mathbb{Z}_p}q^{-y} \lambda ^y [x+y]_q ^n
d\mu _{-q} (y).\label{lamqEp}
\end{eqnarray}
From (\ref{lamqEp}), we can easily derive
$$E_{n,q}(\lambda,x)=\frac{[2]_q}{(1-q)^n}
\underset{l=0}{\overset{n}{\sum}}{n \choose l} (-1)^l
q^{lx}\frac{1}{1+q^l \lambda}.
$$

Let $g(x,\lambda :q)=\underset{n=0}{\overset{\infty}{\sum}}E
_{n,q}(\lambda,x)\frac{t^n}{n!}$. Then we have
\begin{eqnarray*}
g(x,\lambda :q)&=&\underset{n=0}{\overset{\infty}{\sum}}E
_{n,q}(\lambda,x)\frac{t^n}{n!}=\int _ {\mathbb{Z}_p}q^{-y} \lambda
^y  e^{[x+y]_q ^t} d\mu _{-q}
(y)\\
&=&\underset{n=0}{\overset{\infty}{\sum}}\frac{[2]_q}{(1-q)^n}\underset{l=0}{\overset{n}{\sum}}
{n \choose l} (-1)^l q^{lx} (\underset{m=0}{\overset{\infty}{\sum}}
(-1)^m q^{ml}\lambda^m) \frac{t^n}{n!}\\
&=& [2]_q\underset{m=0}{\overset{\infty}{\sum}}(-1)^m \lambda^m
e^{[m+x]_q} t.
\end{eqnarray*}
It follows that
$$E_{n,q} (\lambda,x)=\frac{d^n (g(x, \lambda
:g))}{(dt)^n}|_{t=0}
=[2]_q\underset{m=0}{\overset{\infty}{\sum}}(-1)^m \lambda^m
[m+x]_q^n.$$
Then we obtain the following theorem.

\begin{theorem} For $n \in {\Bbb N} \cup \{ 0 \}$, we have
$$E_{n,q}(\lambda, x)= [2]_q\underset{m=0}{\overset{\infty}{\sum}}(-1)^m \lambda^m
[m+x]_q^n.$$
\end{theorem}

By the same method, we consider the $\lambda$-$q$-Genocchi
polynomials as follows.
\begin{eqnarray}
t {\int_{\Bbb Z_p}}q^{-x} \lambda^x e^{[x+y]_{q} t} d\mu_{-q}(x )
&=&[2]_q t \sum_{m=0}^\infty (-1)^m \lambda^m e^{[m+x]_q t}\label{lamqGp}\\
&=&\sum_{n=0}^\infty G_{n,q}(\lambda, x)\frac{t^n}{n!}.\notag
\end{eqnarray}

By (\ref{lamqGp}), we see
\begin{eqnarray}
E_{n,q}(\lambda, x)= {\int_{\Bbb Z_p}}q^{-y} \lambda^y [x+y]_{q}^n
d\mu_{-q}(y)= \frac{G_{n+1, q}(\lambda, x) }{n+1} \label{EG}
\end{eqnarray}
and $G_{0,q} (\lambda,x)=0$. From the definition of
$\lambda$-$q$-Euler polynomials, we derive
\begin{eqnarray}
E_{n,q}(\lambda, x)&=&  {\int_{\Bbb Z_p}}q^{-y} \lambda^y
[x+y]_{q}^n
d\mu_{-q}(y) \notag\\
&=& \sum_{l=0}^n \binom{n}{l} [x]_q^{n-l} q^{lx} {\int_{\Bbb Z_p}}q^{-y} \lambda^y [y]_{q}^l  d\mu_{-q}(y) \label{Enq}\\
&=&\frac{[2]_q}{[2]_{q^d}} [d]_q^n \sum_{a=0}^{d-1} (-1)^a \lambda^a
{\int_{\Bbb Z_p}}q^{-dy} \lambda^{dy} [\frac{x+a}{d}+y]_{q^d}^n
d\mu_{-q}(y,) \notag
\end{eqnarray}
for $d \in \Bbb{N}$ with $d \equiv 1 \, (\mod 2)$. By (\ref{Enq}),
we see that
\begin{eqnarray*}
E_{n,q}(\lambda, x)=
&=& \sum_{l=0}^n \binom{n}{l} [x]_q^{n-l} q^{lx} E_{l,q} (\lambda, x)\\
&=&\frac{[2]_q}{[2]_{q^d}} [d]_q^n \sum_{a=0}^{d-1} (-1)^a \lambda^a
E_{n, q^d} (\lambda^d, \frac{x+a}{d}). \notag
\end{eqnarray*}

It is easy to see that
$$
qI_{-q} (f_1)+I_{-q} (f)=[2]_q f(0),
$$
where $f_1 (x)=f(x+1)$. Thus, we have
\begin{eqnarray*}
q {\int_{\Bbb Z_p}}q^{-y-1} \lambda^{y+1} [x+1+y]_{q}^n
d\mu_{-q}(y)+{\int_{\Bbb Z_p}}q^{-y} \lambda^{y} [x+y]_{q}^n
d\mu_{-q}(y)=[2]_q [x]_q^n. \notag
\end{eqnarray*}

Therefore, we obtain the following theorem.

\begin{theorem} For $n \in {\Bbb N} \cup \{ 0 \}$, we have
$$\lambda E_{n,q}(\lambda, x+1)+E_{n,q}(\lambda, x)= [2]_q [x]_q^n.$$
\end{theorem}
By Theorem 10 and (\ref{EG}), we have the following result.
\begin{corollary}
For $n \in {\Bbb N} \cup \{ 0 \}$, we have
$$\lambda G_{n,q}(\lambda, x+1)+G_{n,q}(\lambda, x)= [2]_q n [x]_q^{n-1}.$$
\end{corollary}

It is easy to see that
\begin{eqnarray*}
\frac{\partial}{\partial x} [x+y]_q^n &=& n[x+y]_q^{n-1} \frac{\log q}{q-1} q^{x+y}   \\
&=& n \log q \,[x+y]_q^{n-1}+ \frac{\log q}{q-1} n [x+y]_q^{n-1}.
\end{eqnarray*}
From (\ref{lamqEp}), we note that
\begin{eqnarray}
\frac{\partial}{\partial x} E_{n,q}(\lambda,x)=\frac{\partial}{\partial x} {\int_{\Bbb Z_p}}q^{-y} \lambda^y
[x+y]_{q}^n d\mu_{-q}(y). \label{partx}
\end{eqnarray}
The right side of (\ref{partx}) equals
\begin{eqnarray*}
n \log q {\int_{\Bbb Z_p}}q^{-y} \lambda^y [x+y]_{q}^n d\mu_{-q}(y)
+ \frac{\log q}{q-1} n {\int_{\Bbb Z_p}}q^{-y} \lambda^y [x+y]_{q}^n
d\mu_{-q}(y)\\ = n \log q \, E_{n,q}(\lambda, x) + \frac{\log
q}{q-1}\, n E_{n-1,q}(\lambda, x).
\end{eqnarray*}
Therefore, we obtain the following lemma.

\begin{lemma} For $n \in {\Bbb N}$, we have
$$\frac{\partial}{\partial x} E_{n,q}(\lambda,x)=n \log q \, E_{n,q}(\lambda, x) + \frac{\log q}{q-1}\, n
E_{n-1,q}(\lambda, x).$$
\end{lemma}

\vskip 10pt

\begin{remark} Note that
\begin{eqnarray*} \frac{\partial}{\partial x}
G_{n,q}(\lambda, x)&=& nE_{n-1, q}
(\lambda,x) \\
&=& \frac{n}{(1-q)^{n-1}} [2]_q \sum_{l=0}^{n-1} \binom{n-1}{l}
(-1)^l q^{lx} \frac{1}{1+ q^l \lambda}.
\end{eqnarray*}
\end{remark}

\begin{remark} Note that
\begin{eqnarray*}
E_{n,q}(\lambda, dx)&=&{\int_{\Bbb Z_p}}q^{-y} \lambda^y
[dx+y]_{q}^n d\mu_{-q}(y) \\
&=& [d]_q^n \frac{[2]_q}{\,\,[2]_{q^d}}  \sum_{a=0}^{n-1} (-1)^a
\lambda^a \int_{\Bbb Z_p} [x+\frac{a}{d}+y]_{q^d}^n \lambda^{dy}
q^{-dy} d\mu_{-q^d}(y) \\
&=& [d]_q^n \frac{[2]_q}{\,\,[2]_{q^d}}  \sum_{a=0}^{n-1} (-1)^a
\lambda^a E_{n, q^d} (\lambda^d, x+ \frac{a}{d} ),
\end{eqnarray*}
\end{remark}
for $d \in \Bbb{N}$ with $d \equiv 1(\mod 2)$.

\vskip 20pt

For $n \in \Bbb{N}$, it is known that
\begin{eqnarray}
q^n I_{-q} (f_n) =(-1)^n I_{-q} (f) +[2]_q \sum_{l=0}^{n-1}
(-1)^{n-1-l} q^l f(l), \quad (\text{see  [7]}), \label{qI}
\end{eqnarray}
where $f_n (x)=f(x+n)$. By (\ref{qI}), we obtain the following
lemma.

\begin{lemma} For $n \in {\Bbb N}$, we have
$$q^n I_{-q} (f_n) +(-1)^{n-1} I_{-q} (f) = [2]_q \sum_{l=0}^{n-1}
(-1)^{n-1-l}q^l f(l).$$
\end{lemma}
For $n \equiv 1 (\mod 2)$, we also have
$$q^n I_{-q} (f_n) + I_{-q} (f) = [2]_q \sum_{l=0}^{n-1}
(-1)^{l}q^l f(l).$$

If we take $f(x)=\lambda^x q^{-x} [x]_q^m$ with $m \in \Bbb{N}\cup \{ 0 \}$, then we see that
$$q^n {\int_{\Bbb Z_p}}q^{-x-n} \lambda^{x+n}
[x+n]_{q}^m d\mu_{-q}(x)+ {\int_{\Bbb Z_p}}q^{-x} \lambda^x
[x]_{q}^m d\mu_{-q}(x) = [2]_q \sum_{l=0}^{n-1} (-1)^l\lambda^l
[l]_q^m.$$
Thus we have
$$E_{m,q} (\lambda,n)+E_{m,q} (\lambda)=[2]_q \sum_{l=0}^{n-1} (-1)^l \lambda^l [l]_q^m .
$$
For $m \equiv 1 (\mod 2)$, we note that
\begin{eqnarray*}
E_{n,q}(\lambda, dx)
&=& \frac{[2]_q}{\, \, \, [2]_{q^m}} [m]_q^n
\sum_{a=0}^{m-1} (-1)^a \lambda^a
E_{n, q^m} (\lambda^m, \frac{a}{m} )\\
&=& \frac{[2]_q}{\, \, \,[2]_{q^m}} [m]_q^n  \sum_{a=0}^{n} (-1)^a \lambda^a
\sum_{l=0}^{n} \binom{n}{l} (\frac{[a]_q}{[m]_q})^{n-l} q^{la} E_{n, q^m} (\lambda^m,
\frac{a}{m} )\\  &=& \frac{[2]_q}{\, \, \,[2]_{q^m}} \sum_{l=0}^{n}  \binom{n}{l} [m]_q^l E_{l, q^m} (\lambda^m)
\sum_{a=0}^{m-1} (-1)^a \lambda^a q^{la} [a]_q^{n-l} .
\end{eqnarray*}

\begin{remark} Note that
$$
\frac{G_{m+1, q}(\lambda, n) }{m+1} + \frac{G_{m+1, q}(\lambda) }{m+1}
=[2]_q  \sum_{l=0}^{n-1} (-1)^a \lambda^a  [l]_q^m .
$$
\end{remark}

\vskip 10pt

Now we can also consider the following DC type $\lambda$-$q$-Euler numbers and polynomials.
For $\lambda \in \Bbb{C}_p$ with $|1-\lambda|_p<1$, we define the DC type $\lambda$-$q$-Euler numbers as
\begin{eqnarray*}
E_{n,q}^*(\lambda) &=& {\int_{\Bbb Z_p}}\lambda^x [x]_q^n  d\mu_{-q}(x)\\
&=&  \frac{[2]_q}{(1-q)^n} \sum_{l=0}^{n}  \binom{n}{l}  (-1)^l
\frac{1}{1+q^{l+1}\lambda}\\
&=& [2]_q \sum_{m=0}^{\infty} (-1)^m \lambda^m q^m [m]_q^n.
\end{eqnarray*}

Let $g^* (t, \lambda:q)=  \sum_{n=0}^{\infty} E_{n,q}^*(\lambda) \frac{t^n}{n!}$. Then we see that
$$
g^* (t, \lambda:q)={\int_{\Bbb Z_p}}\lambda^x e^{[x]_q t}  d\mu_{-q}(x)
=[2]_q  \sum_{m=0}^{\infty} (-1)^m \lambda^m q^m  e^{[m]_q t}.
$$

The DC type $\lambda$-$q$-Euler polynomials are also defined as
\begin{eqnarray*}
E_{n,q}^*(\lambda,x) &=& {\int_{\Bbb Z_p}}\lambda^y [x+y]_q^n  d\mu_{-q}(y)\\
&=&  \frac{[2]_q}{(1-q)^n} \sum_{l=0}^{n}  \binom{n}{l}  (-1)^l
q^{lx} \frac{1}{1+q^{l+1}\lambda}\\
&=& [2]_q \sum_{m=0}^{\infty} (-1)^m \lambda^m q^m [m+x]_q^n.
\end{eqnarray*}

Thus we can give the generating function of the DC type $\lambda$-$q$-Euler polynomials as follows.
$$
 \sum_{n=0}^{\infty} E_{n,q}^*(\lambda,x)\frac{t^n}{n!}= {\int_{\Bbb Z_p}}\lambda^y e^{[x+y]_q t}  d\mu_{-q}(y)
=[2]_q  \sum_{m=0}^{\infty} (-1)^m \lambda^m q^m  e^{[m+x]_q t}.
$$

\bigskip

\section{Further Remarks and Observation for the $q$-extension Lerch zeta function}
\vskip 10pt

In this section, we assume that $q\in \Bbb{C}$ with $|q|<1$. It is well-known that Lerch type zeta function is defined as
$$\zeta(x,s,a)=\sum_{n=0}^{\infty} \frac{x^n}{(n+a)^s},   $$
where $a \in \Bbb{C}$ with $a \neq 0,-1,-2,\, \cdots$, and $s \in \Bbb{C}$ when $|x|<1$, $Re(S)>1$ when $|x|=1$, and Hurwitz zeta
function is defined as
$$\zeta(s,a)=\sum_{n=0}^{\infty} \frac{1}{(n+a)^s},   $$
where $Re(S)>1$ and $a \neq 0,-1,-2,\, \cdots$.
The Lerch zeta functin is known that
$$\zeta(s,\eta)=\sum_{n=0}^{\infty} \frac{e^{2 \pi i \eta}}{n^s}=e^{2 \pi i \eta} \zeta(e^{2 \pi i \eta},s,1),   $$
where $\eta \in \Bbb{R}$ and $Re(S)>1$.

Now we consider the first kind of the $q$-extension of Lerch type zeta function as follows.
\begin{eqnarray}
\zeta_q (\lambda,s)=(1-q) \frac{2-s}{s-1} \sum_{m=1}^{\infty}
\frac{q^m \lambda^m}{\, [m]_q^{s-1}}+\sum_{m=1}^{\infty} \frac{q^m \lambda^m}{[m]_q^{s}},
\label{lerch}
\end{eqnarray}
where  $q \in \Bbb{C}$ with $|q|<1$, and $\lambda \in \Bbb{C}$ with $\lambda=e^{2 \pi i /f}$ $(f \in \Bbb{N})$.

By Theorem 3, we see that
\begin{eqnarray}
- \frac{\beta_{k,q} (\lambda)}{k}=(q-1) \frac{k+1}{k}
\sum_{m=1}^{\infty} q^m \lambda^m [m]_q^k + \sum_{m=1}^{\infty}  q^m
\lambda^m [m]_q^{k-1},  \label{beta}
\end{eqnarray}
for $k \in \Bbb{N}$.

By (\ref{lerch}) and (\ref{beta}), we obtain the following theorem.

\begin{theorem} For $k \in {\Bbb N}$, we have
$$\zeta_q (\lambda,1-k)=- \frac{\beta_{k,q} (\lambda)}{k}. $$
\end{theorem}

Now, we define the second of the $q$-extension of Lerch zeta function as follows.
For $s \in \Bbb{C}$ and $\lambda=e^{2 \pi i /f}$ $(f \in \Bbb{N})$,
define
\begin{eqnarray}
\zeta_q^* (\lambda,s)=\sum_{m=1}^{\infty} \frac{q^m  \lambda^m}{[m]_q^s}.   \label{30}
\end{eqnarray}
By Theorem 6, we easily see that
\begin{eqnarray}
- \frac{\beta_{k,q} (\lambda)}{k}= \sum_{m=1}^{\infty} q^m \lambda^m
[m]_q^{k-1}.  \label{31}
\end{eqnarray}

By (\ref{30}) and (\ref{31}), we obtain the following theorem.
\begin{theorem} For $k \in {\Bbb N}$, we have
$$\zeta_q^* (\lambda,1-k)=- \frac{\beta_{k,q} (\lambda)}{k}. $$
\end{theorem}

\vskip 10pt

\begin{remark} The extension of Hurwitz's type $q$-Euler zeta
function is defined as
$$
\zeta_{q,E} (\lambda,s)=[2]_q \sum_{m=0}^{\infty} \frac{(-1)^m  \lambda^m}{[m+x]_q^s},
$$
where $s \in \Bbb{C}$, $\lambda \in \Bbb{C}$ with $\lambda=e^{2 \pi i /f}$ $(f \in \Bbb{N})$.
Then we have
$$\zeta_q (\lambda,1-k)=E_{k,q} (\lambda, x), \quad (k \in \Bbb{N}). $$
\end{remark}

\end{document}